\documentclass[12pt]{elsarticle}

\usepackage{lineno,hyperref}
\modulolinenumbers[5]
\usepackage{diagbox}
\usepackage{amsmath,amssymb}

\usepackage[T1]{fontenc}
\usepackage{lmodern}
\usepackage[T1]{fontenc}
\usepackage{lineno,hyperref}
\modulolinenumbers[5]
\usepackage{enumitem}
\usepackage{array}
\usepackage{tabularx} 
\usepackage{adjustbox} 
\usepackage{makecell}
\usepackage{amsmath}
\usepackage{setspace}
\usepackage{float}
\usepackage{graphicx}
\usepackage{subfigure}
\usepackage[utf8]{inputenc}
\usepackage{mathabx}
\usepackage{bm}
\usepackage{mathptmx}
\newsavebox{\foobox}
\newcommand{\slantbox}[2][0]{\mbox{%
        \sbox{\foobox}{#2}%
        \hskip\wd\foobox
        \pdfsave
        \pdfsetmatrix{1 0 #1 1}%
        \llap{\usebox{\foobox}}%
        \pdfrestore
}}
\newcommand\unslant[2][-.18]{\slantbox[#1]{$#2$}}

\journal{Journal of \LaTeX\ Templates}
\usepackage{makecell}
\usepackage{geometry}
\begin{document}

\begin{frontmatter}

\title{\Large \textbf{A Stochastic Multivariate Latent Variable Model For Categorical Responses}}

\author{Mahdi Mollakazemiha$^{\ast}$}

\cortext[mycorrespondingauthor]{Corresponding author}
\ead{mollakazemiha.mahdi@gmail.com}

\address{Department of Statistics, Shahid Beheshti University, Tehran, Iran}

\begin{large}
\begin{abstract}
\setstretch{1.15}

This paper introduces a mathematical framework of a stochastic process model as a generalization of diffusion stochastic processes to model latent variables in categorical responses given unobserved random effects and maximum likelihood estimation of parameters is indicated. 

\end{abstract}
\end{large}


\end{frontmatter}

\begin{large}
\section{Model and Likelihood}
\onehalfspacing

We use $X^{(k)}(t)=(X_1^{(k)}(t),...,X_p^{(k)}(t))'$ to denote the $p-$ dimensional multivariate  Brownian motion process for the $k~$th individual, $k=1,...,m$, with drift and diffusion parameters given by $(\mu^{(k)},\Sigma)$ and expressed with a stochastic differential equation:
\begin{equation}
dX^{(k)}(t)=\mu^{(k)} dt+\Sigma^{(k)\frac{1}{2}} dW^{(k)}(t)
\end{equation}
Where  $W^{(k)}$ is a standard $p-$ dimensional multivariate Brownian motion
process $W^{(k)}$ and $X^{(k)}(t)$ is a multivariate normal distribution
with $\mu^{(k)}=(\mu_1^{(k)},...,\mu_p^{(k)}) \in R^{p }$   mean and $\Sigma$ is
${p\times p}$ covariance matrix. Component of $X(t)$ at times
$t_0<t_1<t_2<...<t_n$  is observed. For the moment assume that the time points $t_0$, ..., $t_n$ are equally spaced, and
denote the common time increment as $(t_i-t_{i-1})$. 
The natural approach taken by many authors is to use an Euler approximation of the process which is of the form
$$X^{(k)}(t_i-t_{i-1})|X^{(k)}(t_0) \sim N (X^{(k)}(t_0) +(t_i-t_{i-1})\mu_j^{(k)}, (t_i-t_{i-1})\Sigma),$$
where $X^{(k)}(t_0)=0$.

We use $X^{(k)}_{ij}=X^{(k)}_{j}(t_i)$ to denote $j~$th the thought process outcome for the $k~$th individual at time $t_i$ and let $Y^{(k)}_{ij}=X^{(k)}_{j}(t_i)-X^{(k)}_{j}(t_{i-1})$, $k=1,...,m$, $i=1,...,n$ and $j=1,..., p$.  The model for thought process takes the form:
\begin{equation}
    Y^{(k)}_{ij}=(t_i-t_{i-1})\mu^{(k)}_{j}+\varepsilon_{ijk}^{(1)},~~~~i=1,...,n~,~~j=1,..., p,
\end{equation}
Where $\mu^{(k)}_{ij}={V_{jk}^{(1)}}'\beta_j^{(1)}+{U_{jk}^{(1)}}'b_k^{(1)}$. Also,
$\varepsilon_{ik}^{(1)}=(\varepsilon_{i1k}^{(1)},...,\varepsilon_{ipk}^{(1)})'$ is  multivariate normal  with zero mean and $\Sigma$ covariance matrix and  $V_{jk}^{(1)}$ denote the  vector of  explanatory covariates for $Y_{ij}$ and and $U^{(1)}_{j}$ be some sub-vectors of covariate of $V^{(1)}_{jk}$. Also, $b_k^{(l)}$ is the vector of random effect where
$$\begin{array}{l} 
{ b_k^{(1)}\mathop \sim 
\limits^{iid} MVNormal(0},\Sigma_1),\\  
\end{array}$$ 
 To our model, one may 
use random effects for outcomse across time, which leads 
to the conditional independence of the vector of responses in 
different occasions given subject-level effects $b_k^{(1)}$. 

Let $R^k_{j}$ be the reaction time for $j~$th the thought process outcome for the $k~$th individual  and $V_{jk}^{(2)}$ denote the vector of  explanatory covariates for $R^k_{j}$. Suppose that we are interested in making inferences about the effect of $V_{jk}^{(2)}$ on the outcome  of the reaction time. 
One implication of this model is on the log scale, directly in terms of $V_{jk}^{(2)}$ and commonly be written as

\begin{equation}
 log(R^k_{j}) ~ | \{ Y^{(k)}_{ij}< a_1 ~ or ~ Y^{(k)}_{ij}> a_2\}={V_{jk}^{(2)}}'\beta_j^{(2)}+{U_{jk}^{(2)}}'b_k^{(2)}+\varepsilon_{ijk}^{(2)},
\end{equation}
where $a_1$ and $a_2$ are the parameters representing boundaries of decision making in the model and $U^{(2)}_{jk}$ be some sub-vectors of covariate of $V^{(2)}_{jk}$. Also, $\varepsilon_{ik}^{(2)}=(\varepsilon_{i1k}^{(2)},...,\varepsilon_{ipk}^{(2)})'$ is  multivariate normal  with zero mean and $\Sigma_{R}$ covariance matrix and $b_k^{(2)}$ is the vector of random effect where
$$\begin{array}{l} 
{ b_k^{(2)}\mathop \sim 
\limits^{iid} MVNormal(0},\Sigma_2),\\  
\end{array}$$ 
\subsection{The Joint Model}
The joint model is assumed to take the form:
$$\begin{array}{l}
 Y^{(k)}_{ij}=(t_i-t_{i-1})\mu^{(k)}_{j}+\varepsilon_{ijk}^{(1)},~~~~i=1,...,n~,~~j=1,..., p,~k=1,..., m\\
 log(R^k_{j}) ~ |\{ ~  Y^{(k)}_{ij}< a_1 ~ or ~ Y^{(k)}_{ij}> a_2\}={V_{jk}^{(2)}}'\beta_j^{(2)}+{U_{jk}^{(2)}}'b_k^{(2)}+\varepsilon_{ijk}^{(2)},
\end{array}$$
Where $\mu^{(k)}_{ij}={V_{jk}^{(1)}}'\beta_j^{(1)}+{U_{jk}^{(1)}}'b_k^{(1)}$ and
$$\begin{array}{l}
{( b_k^{(1)},b_k^{(2)})\mathop  \sim
\limits^{iid} MVNormal(0},\Sigma_{B}),\\
\Sigma_{B}= \left( {\begin{array}{*{20}c}
   {\Sigma _{1} } & {\Sigma _{12} }  \\
   {\Sigma _{21} } & {\Sigma _{2} }  \\
\end{array}} \right),
\end{array}$$
where $\Sigma_{1}=Var(b_k^{(1)})$,  $\Sigma_{2}=Var(b_k^{(2)})$ and  $\Sigma_{12}=\Sigma_{21}'=Cov(b_k^{(1)},b_k^{(2)})$. 

To joint model the  bivariate vector of  outcomes, one may
 use random effects to take into account the
correlation between outcomes across time, which leads
to the conditional independence of the vector of outcomes in
different occasions given  subject-level effects $b_k^{(1)}$ and $b_k^{(2)}$.

Because
of identifiability problem we have to assume
 $\Sigma=\Sigma_R=I$, ($I$ is a diagonal matrix).
The vector of parameters $\beta_{j}^{(1)}$, $\beta_{j}^{(2)}$, the parameters of  $a_1$, $a_2$ and $\Sigma_B$ 
should be estimated. Maximum likelihood estimates for the parameters can be obtained with commonly used algorithms for maximizing the likelihood.

\subsection{Likelihood}
The log-likelihood function under the joint model  is
\begin{equation}
    \begin{array}{l} 
    log ~L=\sum_{k = 1}^m log[{\int\limits_{b_k^{(1)} } \int\limits_{b_k^{(2)} }{\prod\limits_{j = 
    1}^{p } \prod\limits_{i = 1}^{n} {\mathbb{P}(log(R_{j}^{k})=r_{ij}^{k}| Y_{ij}^{(k)}< a_1~  or ~ Y_{ij}^{(k)}> a_2, b_k^{(2)})}}} \\
    ~~~~~~\mathbb{P}(Y_{ij}^{(k)}< a_1~  or ~ Y_{ij}^{(k)}> a_2| {b_k^{(1)} }))~\varphi_{12} (b_k^{(1)} ,b_k^{(2)} ;\Sigma_1,\Sigma_2)\;db_k^{(1)}db_k^{(2)}],
    \end{array}
\end{equation}
where  $\mathbb{P}(R_{j}^{k}=r_{j}^{k} ~ \vert ~ Y_{ij}^{(k)}<a_{1} ~ or ~ Y_{ij}^{(k)}>a_{2}, b^{(2)}_{k})$ denotes the corresponding density functions for $R_{j}^{k}$ given $\{Y_{ij}^{(k)}< a_{1} ~ or ~ Y_{ij}^{(k)} > a_{2}\}$ and $b_k^{(2)}$. We have:
\begin{equation}
    \begin{array}{l} 
    \mathbb{P}(R_{j}^{k} = r_{j}^{k} ~ \vert ~ Y_{ij}^{(k)} < a_1~  or ~ Y_{ij}^{(k)}> a_2, b_k^{(2)}) \\
    = \frac{\partial}{\partial{r_{j}^{k}}} \mathbb{P}(R_{j}^{k} \le r_{j}^{k} ~ \vert ~ Y_{ij}^{(k)}< a_1~  or ~ Y_{ij}^{(k)}> a_2, b_k^{(2)})\\
    = \frac{\partial}{\partial{r_{j}^{*k}}} \mathbb{P}(log(R_{j}^{k}) \le log(r_{j}^{*k}) ~ \vert ~ Y_{ij}^{(k)}< a_1~  or ~ Y_{ij}^{(k)}> a_2, b_k^{(2)})
    \end{array}
\end{equation}

Where $r^{*}_{j} = log(r_{j})$ and $\frac{\partial}{\partial{r_{j}^{*k}}}(.)$ is derived from $r_{j}^{*k}$. Also, $\phi_{12}(.;.)$ denotes the corresponding joint density for $b^{(1)}_{k}$ and $b^{(2)}_{k}$.
Also, we have $\mathbb{P}(Y_{ij}^{(k)}<a_{1}) = \Phi (m_{1}) $ ~and~ $\mathbb{P}(Y_{ij}^{(k)}>a_{2}) = (1-\Phi(m_{2}))$, $m_{1}= \frac{a_{1}-({V_{jk}^{(1)}}'\beta_j^{(1)}+{U_{jk}^{(1)}}'b_k^{(1)})}{Var(Y_{ij}^{k})}$,~~ $m_{2} = \frac{a_2-({V_{jk}^{(1)}}'\beta_j^{(1)}+{U_{jk}^{(1)}}'b_k^{(1)})}{Var(Y_{ij}^{k})}$.

To obtain $\mathbb{P}(R_{j}^{k} = r_{j}^{k}  ~ \vert  ~ Y_{ij}^{k}< a_{1}~ or~ Y_{ij}^{k}>a_{2}, b^{(2)}_{k} )$, we use an approximation based on conditional expectations and regression with binary variables. (Joe, 1995), proposed new approximations for multivariate normal probabilities for rectangular regions based on conditional expectations and regression with binary variables. 

Let $( Y_{ij}^{(k)},log R_{j}^{k})$ given $(b^{(1)}_{k}, b^{(2)}_{k})$ be bivariate normal random vector with mean vector equals to  ${V_{jk}^{(1)}}'\beta_j^{(1)}+{U_{jk}^{(1)}}'b_k^{(1)}$ and ${V_{jk}^{(2)}}'\beta_j^{(2)}+{U_{jk}^{(2)}}'b_k^{(2)}$ with matrix covariate $\Sigma$ and $\Sigma_R$. Also, Considering ~ \( \mathbb{I}(R_{j}^{k}\le r_{j}^{*k}),  ~ \mathbb{I}(Y_{ij}^{(k)} < a_{1}), ~ \mathbb{I}(Y_{ij}^{(k)} > a_{2}) \) are indicator functions, we have:

\begin{equation}
    \begin{array}{l}
    \mathbb{P}(log(R_{j}^{k})=r_{j}^{k}| Y_{ij}^{(k(}< a_1~  or ~ Y_{ij}^{(k(}> a_2, b_k^{(2)}) \\
    = \mathbb{E}[ ~ \mathbb{I}(log(R_{j}^{k})\le r_{j}^{*k}) ~ \vert ~ \mathbb{I}_{(Y_{ij}^{(k)}< a_1)} =1 ~  or ~ \mathbb{I}_{(Y_{ij}^{(k)}> a_2)} = 1 ~,~ b_{k}^{2} ]
    \end{array}
\end{equation}

Based on Joe's approximation we get:

\begin{equation}
    \begin{array}{l}
    \mathbb{E}[ ~ \mathbb{I}(log(R_{j}^{k})\le r_{j}^{*k}) ~ \vert ~ \mathbb{I}_{(Y_{ij}^{(k)}< a_1)} =1 ~  or ~ \mathbb{I}_{(Y_{ij}^{(k)}> a_2)} = 1 ~,~ b_{k}^{2} ] \\ 
    \cong ~ \mathbb{E}[~\mathbb{I}(~log(R_{j}^{k})\le r_{j}^{*k}~)~ ] + \Omega_{21}\Omega^{-1}_{11}(1- \mathbb{I}_{(Y_{ij}^{(k)}< a_1) ~ \cup ~  (Y_{ij}^{(k)}> a_2)})^{T}
    \end{array}
\end{equation}

where  \( \Omega_{21} = Cov( ~ \mathbb{I}_{log(R_{j}^{k})} \le r^{k*}_{j}, ~ \mathbb{I}_{(Y_{ij}^{(k)}< a_1) ~ \cup ~  (Y_{ij}^{(k)}> a_2)} ~ ) \) ~ and ~  \( \Omega_{11} = Var( ~ \mathbb{I}_{(Y_{ij}^{(k)}< a_1) ~ \cup ~  (Y_{ij}^{(k)}> a_2)} ~ ) \), and \(\mathbb{E}[~\mathbb{I}(~log(R_{j}^{k})\le r^{k*}_{j}~)~ ]  = \mathbb{P}(log(R_{j}^{k}) < r^{k*}_{j}) \). \\

\(\Omega_{21}\) can be written as following:

\begin{equation}
    \begin{array}{l} 
       \Omega_{21}  = \mathbb{E}[~ \mathbb{I}_{(log(R_{j}^{k}) \le r^{k*}_{j})} ~ . ~ \mathbb{I}_{(Y_{ij}^{(k)}< a_1) ~ \cup ~  (Y_{ij}^{(k)}> a_2)} ~] \\
        - \mathbb{E}[~  \mathbb{I}_{(log(R_{j}^{k}) \le r^{k*}_{j})}~] ~ . ~  \mathbb{E}[~ \mathbb{I}_{(Y_{ij}^{(k)}< a_1) ~ \cup ~  (Y_{ij}^{(k)}> a_2)} ~] \\
         = \mathbb{E}[~ \mathbb{I}_{(log(R_{j}^{k}) \le r^{k*}_{j}) ~ \cap ~ [ (Y_{ij}^{(k)}< a_1) ~ \cup ~  (Y_{ij}^{(k)}> a_2)] } ~] \\
         - \mathbb{P}(log(R_{j}^{k}) \le r^{k*}_{j})[~ \mathbb{P}(Y_{ij}^{(k)} < a_{1}) + \mathbb{P}(Y_{ij}^{(k)} > a_{2}) ~] \\
         = \mathbb{P}(log(R_{j}^{k}) \le r^{k*}_{j} ~,~  Y_{ij}^{(k)} < a_{1} ~or~ Y_{ij}^{(k)} > a_{2}~) \\
         - \mathbb{P}(log(R_{j}^{k}) \le r^{k*}_{j})[~ \mathbb{P}(Y_{ij}^{(k)} < a_{1}) + \mathbb{P}(Y_{ij}^{(k)} > a_{2}) ~] 
    \end{array}
\end{equation}

and for \(\Omega_{11}\) we have:

\begin{equation}
    \begin{array}{l} 
       \Omega_{11}  = \mathbb{E}[~ \mathbb{I}^{2}_{(Y_{ij}^{(k)}< a_{1}) ~\cup~ (Y_{ij}^{(k)} > a_{2})}  ~] - \mathbb{E}^{2}[~ \mathbb{I}_{(Y_{ij}^{(k)}<a_{1} ~\cup~ Y_{ij}^{(k)}>a_{2})} ~] \\
        = \mathbb{E}[~ \mathbb{I}_{(Y_{ij}^{(k)}<a_{1}) ~\cup~ (Y_{ij}^{(k)}>a_{2})} ~] - \mathbb{E}^{2}[~ \mathbb{I}_{(Y_{ij}^{(k)}<a_{1} ~\cup~ Y_{ij}^{(k)}>a_{2})} ~] \\
       = \mathbb{P}(Y_{ij}^{(k)}<a_{1}) + \mathbb{P}(Y_{ij}^{(k)}>a_{2}) \\
       - \mathbb{P}(Y_{ij}^{(k)}<a_{1})^{2} - \mathbb{P}(Y_{ij}^{(k)}>a_{2})^{2} - 2\mathbb{P}(Y_{ij}^{(k)}<a_{1})\mathbb{P}(Y_{ij}^{(k)}>a_{2}),
    \end{array}
\end{equation}
\newcommand{\comment}[1]{}
\comment{ 
\begin{align*}
   f_{R_{j}}(r_{j} ~ \vert ~ \unslant\theta) + \{~f_{R_{j},Y_{ij}^{(k)}}(r_{j}, y_{ij}^{(k)} ~ \vert ~ {\unslant\theta}) ~ - ~ f_{R_{j}}(r_{j} ~ \vert ~ \unslant\theta) [\Phi(a_{1}) - \Phi(a_{2})+1] ~\} &\\
   \times~ \{\Phi(a_{1}) - \Phi(a_{2}) +1 - (\Phi(a_{1}) - \Phi(a_{2}) +1)^{2} \}
\end{align*} }
Maximum likelihood estimation as a standard method of parameter estimation is applied. The log-likelihood function maximizes at the points which we are interested in. Such methods commonly are not yield analytically, thus numerical approach should be applied while maximizing the likelihood function. Nelder-Mead Simplex method (Nelder \& Mead, 1965) was conducted for this purpose in this study.

\end{large}

\end{document}